\documentclass[journal]{IEEEtran}

\usepackage[dvips]{graphicx}
\usepackage{amssymb,amsfonts,amsmath}
\usepackage{eucal,bm}
\usepackage{gensymb,subfigure}
\usepackage[square,sort,comma,numbers]{natbib}
\usepackage{booktabs}
\usepackage{times,color,url}
\usepackage{xparse}
\ExplSyntaxOn
\NewDocumentCommand{\mref}{m}{\quinn_mref:n {#1}}
\seq_new:N \l_quinn_mref_seq
\cs_new:Npn \quinn_mref:n #1
 {
  \seq_set_split:Nnn \l_quinn_mref_seq { , } { #1 }
  \seq_pop_right:NN \l_quinn_mref_seq \l_tmpa_tl
  ( % print the left parenthesis
  \seq_map_inline:Nn \l_quinn_mref_seq
    { \ref{##1},\nobreakspace } % print the first references
  \exp_args:NV \ref \l_tmpa_tl % print the last or only one
  ) % print the right parenthesis
 }
\ExplSyntaxOff

\graphicspath{{referee_figures/}}
\begin{document}

\title{Optimal Distributed Control of Reactive Power via the Alternating Direction Method of Multipliers}
%
% \author{
% \authorblockN{Petr \v{S}ulc}
% \authorblockA{R. Peierls Centre for Theoretical Physics \\ University of Oxford \\ Oxford OX1 3NP,  UK
%  \\Email: p.sulc1@physics.ox.ac.uk}
% \and
% \authorblockN{Scott Backhaus}
% \authorblockA{Materials, Physics \& Applications Division\\
% LANL and New Mexico Consortium,\\ Los Alamos, NM
% 87545(4), USA\\
% Email: backhaus@lanl.gov}
% \and
% \authorblockN{Michael Chertkov}
% \authorblockA{CNLS \& Theoretical Divison\\
% LANL and New Mexico Consortium,\\ Los Alamos, NM
% 87545(4), USA\\
% Email: chertkov@lanl.gov}
% }

\author{Petr \v{S}ulc,
       Scott Backhaus,
       Michael Chertkov,~\IEEEmembership{Member,~IEEE}% <-this % stops a space
\thanks{P. \v{S}ulc is with Rudolf Peierls Centre for Theoretical Physics, University of Oxford, Oxford, OX1 3NP, UK, e-mail: p.sulc1@physics.ox.ac.uk}% <-this % stops a space
\thanks{S. Backhaus is with Materials, Physics \& Applications Division, LANL, Los Alamos, NM 87545, USA, e-mail: backhaus@lanl.gov}% <-this % stops a space
\thanks{M. Chertkov is with CNLS and Theoretical Divison, LANL, Los Alamos, NM 87545, USA, e-mail: chertkov@lanl.gov}
\thanks{M. Chertkov and S. Backhaus are also with New Mexico Consortium, Los Alamos, NM 87544, USA}

}

\date{\today}
 \maketitle

\begin{abstract}
We formulate the control of reactive power generation by photovoltaic inverters in a power distribution circuit as a constrained optimization that aims to minimize reactive power losses subject to finite inverter capacity and upper and lower voltage limits at all nodes in the circuit. When voltage variations along the circuit are small and losses of both real and reactive powers are small compared to the respective flows, the resulting optimization problem is convex.  Moreover, the cost function is separable enabling a distributed, on-line implementation with node-local computations using only local measurements augmented with limited information from the neighboring nodes communicated over cyber channels.  Such an approach lies between the fully centralized  and local policy approaches previously considered.  We explore protocols based on the dual ascent method and on the Alternating Direction Method of Multipliers (ADMM) and find that the ADMM protocol performs significantly better.
 \end{abstract}


\section{Introduction}

The generation or consumption of reactive power by inverters has been explored by several researchers as a way to control voltage fluctuations in distribution circuits with a high penetration of distributed photovoltaic (PV) generation.  See \cite{KostyaMishaPetrScott3} and references therein for an overview.  These approaches have tended to fall near two extremes.  One extreme considers centralized optimization \cite{KostyaMishaPetrScott} where computations are done by a single central authority which is assumed to have full observability of the system.  Such a system requires two-way communications between the central authority and at least all of controlled inverters.  The rate of these communications should be sufficient to respond to the fastest expected fluctuations of solar irradiance.  At the other extreme are local policy-based methods that require no communications at all, except for perhaps between devices at a single node.  These inherently suboptimal methods rely only on node-local measurements as inputs to a policy that converts the measurements into a control action.  Such policies have been based on heuristics and physical reasoning \cite{KostyaMishaPetrScott2,KostyaMishaPetrScott3}  and on Monte Carlo-like approaches \cite{kundu2013distributed} that use centralized optimization in off-line computations to find  strong correlations between the local measurements and the optimal local control actions.

In \cite{kundu2013distributed}, it was found that a significant source of the suboptimal behavior of local policy-based methods was saturation of inverters.  In this case, the local policy results in a {\it desired} control action that, when combined with the local power flow conditions, is beyond the capability of the inverter.  In \cite{kundu2013distributed}, this generally occurred when the real power injection by the inverter was approaching the apparent power capacity leaving little room for reactive power generation or consumption.  Instead of achieving the desired reactive power injection, the inverter saturated at its apparent power capacity.  In a centralized approach, optimality is restored because the full observability of the central authority makes it aware of the saturation allowing it to compensate with extra reactive power generation (or consumption) from nearby nodes. Such a response suggests that nearest-neighbor communications may be used to restore optimality in capacity-limited or otherwise constrained systems.  There are other reasons to expect that limited communications may provide significant advantages over policy-based control without the larger overhead of centralized communications.  Primary among these is the ability to adapt to system configurations that were unforeseen in the development of a policy based either on heuristics or sampling methods.

In this manuscript, we explore distributed approaches to control reactive power from PV inverters in distribution circuits and show how to restore optimality and adaptability through an iterative message-passing algorithm.  Although we invoke limited communications, we continue to rely solely on local computations which depend only on local measurements and the most current data communicated from the nearest-neighbor nodes. The suggested cyber-physical control scheme is decentralized but optimal.

The results in this manuscript are based on the observation that the convex optimization formulation of \cite{KostyaMishaPetrScott,KostyaMishaPetrScott2,KostyaMishaPetrScott3} is separable in the key optimization variables, i.e. the node voltages and the power flows along the circuit. Separability suggests application of modern methods of distributed computations such as the dual ascent method and the Alternating Direction Method of Multipliers (ADMM) \cite{Boyd2011,boyd2013}.  % Both methods are provably convergent, however,
The ADMM algorithm converges significantly faster than the dual ascent algorithm, a property one expects from the general arguments expressed in \cite{Boyd2011}.
Furthermore, the dual ascent algorithm for a non-differentiable dual function does not always converge to the optimal solution, but can converge to its neighborhood \cite{bertsekas1999nonlinear,nedic2009distributed}.
Faster convergence of ADMM is largely confirmed in our experiments conducted over seven different distribution circuit configurations with different numbers of nodes and varying photovoltaic penetration and load profiles. In each experiment, we minimize the total loss of real power while constraining all node voltages to be within nominal operational bounds.

%The current work appears to be the first application of the ADMM algorithm to the voltage-constrained loss-minimization problem.
The dual ascent method has been considered by \cite{Zampieri13a, Zampieri13b, Zampieri13c}, where a slightly more general model of radial power flows was considered.\footnote{References \cite{Zampieri13a, Zampieri13b, Zampieri13c} approximate power flows along the lines assuming that voltage variations along the line is much smaller than the voltage magnitude at the head of the line.  The LinDistFlow approximation of \cite{89BWa}, used in \cite{KostyaMishaPetrScott,KostyaMishaPetrScott2,KostyaMishaPetrScott3} and adopted in this manuscript to model power flows, assumes additionally that losses of real and reactive power anywhere along the feeder are much smaller than respective flows.}
Dual decomposition distributed algorithm with gradient ascent for voltage regulation was also proposed in \cite{zhang2012optimal}, where
it was shown to solve a convex relaxation of power flow equations.
We further note that an optimization algorithm for optimal load control for frequency regulation has been recently proposed \cite{zhao2013powerB}, where it was shown that a frequency-based load control together with the system dynamics and power flows act as a decentralized primal-dual algorithm that solves the global optimization problem.
The dual ascent and ADMM-based algorithm was also used to solve a semi-definite programming relaxation of optimal power flow problem, where communication was carried out between different segments of the distribution network \cite{dall2013distributed}.  The semi-definite programming relaxation of power flow problem for optimization of real and reactive PV generation in a radial network was also used in \cite{Dall2013},  where it yields an exact solution of the original problem in a single-phase radial network.
%{\color{red} Petr I have added here \cite{Dall2013} ... please check if it is appropriate.}
After the initial submission of our manuscript to the arxiv, other methods of solving convex relaxation of power flow equations via ADMM were proposed \cite{Erseghe2014,magnusson2014distributed,peng2014distributed}.

The material in the remainder of this manuscript is organized as follows. Power flows in a distribution circuit and control of inverters as a global optimization are reviewed in Section \ref{sec:DF} and Section \ref{sec:optim}, respectively. Algorithms for distributed control based on nearest neighbor communications are described in Section \ref{sec_ADMM} and \ref{sec_noU}.
%Section \ref{sec_noU} for the case without voltage constraints.  These algorithms are generalized to account for the voltage constraints in Section \ref{sec_withU}.
The algorithms are tested and compared in Section \ref{sec_cases}. Section \ref{sec_disc} presents our conclusions and a brief discussion of the path forward.

\section{Distributed Flow Formulation}
\label{sec:DF}

\begin{figure}[b]
\centering
\includegraphics[width=8.5cm]{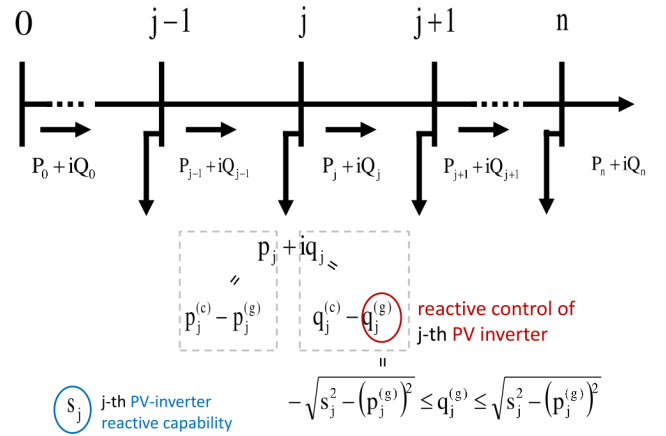}\\~\\
\caption{
A schematic diagram of the distribution circuit illustrating the notation used in Eqs.~(\ref{Pj+1},\ref{Qj+1},\ref{Vj2},\ref{p_q}).
}
\label{fig_feeder}
\end{figure}

The flow of electric power in the quasi-static approximation is governed by Kirchoff's laws. The DistFlow equations \cite{89BWa,89BWb,89BWc} are these laws restated in terms of power flows and applied to radial or tree-like distribution circuit with a discrete set of loads. For the radial case, the DistFlow equations are
\begin{subequations}\label{DistFlow}
$\forall j=0,\ldots,n-1$,
\begin{eqnarray}
&&P_{j+1}\!=\! P_j\!-\!r_j\frac{P_j^2\!+\!Q_j^2}{V_j^2}\!-\!p_{j+1}, \label{Pj+1}\\
&&Q_{j+1}\!=\!Q_j\!-\!x_j\frac{P_j^2\!+\!Q_j^2}{V_j^2}\!-\!q_{j+1}, \label{Qj+1}\\
&&V_{j+1}^2\!=\!V_j^2\!-\!2(r_jP_j\!+\!x_jQ_j)\!+\!(r_j^2\!+\!x_j^2)
\frac{P_j^2\!+\!Q_j^2}{V_j^2},
\label{Vj2}
\end{eqnarray}
where $P_j+iQ_j$ is the complex power flowing away from node $j$ toward node $j+1$, $V_j$ is the voltage at node $j$, $r_j+ix_j$ is the complex impedance of the link between node $j$ and $j+1$, and $p_j+i q_j$ is the complex power extracted at the node $j$. Both $p_j$ and $q_j$ are composed of local consumption minus local generation due to the PV inverter, i.e.
\begin{eqnarray}
p_j=p_j^{(c)}-p_j^{(g)},\quad q_j=q_j^{(c)}-q_j^{(g)}.
\label{p_q}
\end{eqnarray}
Of the four contributions to $p_j+i q_j$, we assume that $p_j^{(g)}$, $p_j^{(c)}$, and $q_j^{(c)}$ are uncontrolled (i.e. driven by consumer load or instantaneous PV generation).  In contrast, the reactive power generated by the PV inverter, $q_j^{(g)}$, can be adjusted within limits.  Eqs.~(\ref{Pj+1},\ref{Qj+1},\ref{Vj2},\ref{p_q}) are solved with the following boundary conditions
\begin{eqnarray}
V_0=\rm{const},\quad P_n=Q_n=0. \label{boundary}
\end{eqnarray}
\end{subequations}
The schematic distribution circuit in Fig.~\ref{fig_feeder} helps to explain the notation.

\section{Control of Inverters as a Global Optimization}
\label{sec:optim}

We aim to solve the following global optimization problem: minimize the total loss of real power while constraining the voltage within nominal operational limits and the reactive power generation to the inverters' apparent power capacity $s_j$:
\begin{subequations}\label{DistFlowOptimization}
\begin{eqnarray}
&\min\limits_{q^{(g)},P,Q,V}& \sum_{j=0}^{n-1} r_j\frac{P_j^2+Q_j^2}{V_j^2},\label{Loss}\\
& \mbox{s.t.}& \mbox{Eq.~(\ref{Pj+1},\ref{Qj+1},\ref{Vj2},\ref{p_q},\ref{boundary})},\nonumber\\
&& \forall j=1,\ldots,n:\quad  \nonumber \\
&& (1-\epsilon)^2 V_0^2 \leq V_j^2\leq (1+\epsilon)^2 V_0^2,\label{Voltage}\\
&& \left| q_j^{(g)} \right| \leq \sqrt{s_j^2- \left(p_j^{(g)} \right)^2}. \label{PV_constraint}
\end{eqnarray}
\end{subequations}
Here, Eq.~(\ref{Voltage}) are the voltage constraints (with $\epsilon$ typically set to $0.05$, following the ANSI C84.1-2006 standard), and Eq.~(\ref{PV_constraint}) is the inverter apparent power constraint.

Under normal operations, the changes in voltage from node to node are small compared to the voltages and the loss of real and reactive power are small compared to the power flows themselves.  In this limit, Eqs.~(\ref{DistFlowOptimization}) can be restated within the LinDistFlow approximation, i.e.,
\begin{subequations}\label{LinDistFlow}
\begin{eqnarray}
&\min\limits_{q^{(g)},P,Q}& \sum_{j=0}^{n-1} r_j\frac{P_j^2 + Q_j^2}{V_0^2}, \label{LossLinear}\\
&\mbox{s.t.}&
\forall j=0,\ldots,n-1: \nonumber \\
&& P_{j+1} = P_j-p_{j+1}^{(c)}+p_{j+1}^{(g)}, \quad P_n=0,\label{P-lin}\\
&& Q_{j+1}=Q_j-q_{j+1}^{(c)}+q_{j+1}^{(g)},  \quad  Q_n=0,\label{Q-lin}\\
&& U_{j+1}=U_j-2(r_jP_j+x_jQ_j),\label{V-lin}\\
&& \forall j=1,\ldots,n: \nonumber \\
&& V_0^2 \left(\epsilon^2 - 2\epsilon \right) \leq U_j\leq V_0^2\left(\epsilon^2 + 2\epsilon\right),\label{Voltage-U}\\
&& \left| q_j^{(g)} \right| \leq \tilde{s}_j, \label{PV_constraint_1}
\end{eqnarray}
where $U_j=V_j^2 - V_0^2$ and $\tilde{s}_j=\sqrt{s_j^2-(p_j^{(g)})^2}$.  We have assumed $V_j^2 \approx V_0^2$ in \eqref{LossLinear}.

Simulations in \cite{KostyaMishaPetrScott,kundu2013distributed} suggest that the LinDistFlow are well justified for a wide range of distribution circuits.  This observation is powerful because the LinDistFlow formulation of Eqs.~(\ref{LinDistFlow}) is convex (a quadratic objective function with linear constraints). Convexity implies that this optimization can be solved efficiently provided each node can communicate with a central authority which performs the computations and distributes the optimal values of $q_j^{(g)}$ to all nodes \cite{KostyaMishaPetrScott,kundu2013distributed}. In the remainder of this work, we will focus on developing a decentralized optimization algorithm, which solves Eqs.~(\ref{LinDistFlow}) only by passing messages between nearest neighbors on the network.

We note that the $P_j$ are determined by solving Eq.~(\ref{P-lin}) for given $p_j^{(c)}$ and $p_j^{(g)}$, and the $P_j$ in Eqs.~(\ref{LossLinear},\ref{V-lin}) can be treated as constants and one can formulate the problem as an optimization over the
$Q_j$ by combining Eqs.~(\ref{Q-lin},\ref{PV_constraint_1}):
\begin{equation}
 \left| Q_j - Q_{j-1} + q_j^{(c)} \right| \leq \tilde{s}_j,  \quad Q_n = 0 \label{PV_constraint_1Q}.
\end{equation}
\end{subequations}
% , i.e.
% \begin{eqnarray}
% \label{LinDistFlowQ}
% &\min\limits_{Q}& \sum_{j=0}^{n-1} \frac{r_j Q_j^2}{V_0^2}, \label{LossLinearQ}\\
% &\mbox{s.t.}& \forall j=1,\cdots,n: \nonumber \\
% && \left| Q_j - Q_{j-1} + q_j^{(c)} \right| \leq \tilde{s}_j,  \quad Q_n = 0 \label{PV_constraint_1Q} \\
% &&  U_{j}=U_{j-1}-2(r_{j-1} P_{j-1}+x_{j-1} Q_{j-1}),\label{V-linQ}\\
% &&  V_0^2 \left(\epsilon^2 - 2\epsilon \right) \leq U_j\leq V_0^2\left(\epsilon^2 + 2\epsilon\right),\label{Voltage-UQ}
% \end{eqnarray}
The actual control outputs are the $q_j^{(g)}$, but these can be inferred from Eq.~(\ref{Q-lin}), once the optimal solution is stated in terms of $Q_j$.

In the following, we present an ADMM-based distributed algorithm for the solution of the LinDistFlow problem.
We further present a dual-ascent distributed algorithm that solves a simplified LinDistFlow problem where the voltage constraints
are omitted and compare its performance with the ADMM algorithm (with and without voltage constraints).

% find it instructive to first address a simplified version of the optimization problem that ignores the voltage constraints, i.e. a problem where we are certain that the voltage resulting from the  optimal $q_j^{(g)}$ will be within the feasible region of Eq.~\ref{Voltage-U}. In this simpler problem, we omit both the LinDistFlow equation (Eq.~\ref{V-lin}) and the voltage constraint (Eq.~\ref{Voltage-U}) yielding:
% \begin{eqnarray}
% \label{LinDistFlowQnoU}
% &\min\limits_{Q}& \sum_{j=0}^{n-1} \frac{r_j Q_j^2}{V_0^2}, \label{LossLinearQnoU} \\
% &\mbox{s.t.}&  \forall j=1,\ldots,n: \nonumber \\
% && Q_j - Q_{j-1} + q_j^{(c)}  - \tilde{s}_j \leq 0,\label{PV_constraint_1A} \\
% && -Q_j + Q_{j-1} - q_j^{(c)}  - \tilde{s}_j \leq 0,  \quad Q_n = 0.  \label{PV_constraint_1B}
% \end{eqnarray}
%
% Two distributed algorithms that solve the  simplified optimization problem \eqref{LossLinearQnoU} are proposed in Section \ref{sec_noU}. Voltage constrains are reintroduced in Section \ref{sec_withU} along with a distributed algorithm for solving the original problem with the voltage constraints included. %\eqref{LossLinearQ}.

\subsection{ADMM consensus distributed algorithm}
\label{sec_ADMM}

We adapt a consensus version of the ADMM algorithm to our problem. A general discussion of the method and proof of convergence is described in \cite{Boyd2011}. The consensus version assumes that each of the nodes in the network has its own local objective function and a local set of constraints which act on a global variable shared between all the nodes. Each node solves a local optimization problem for respective local copies of the global variables. The local optimization problem consists in finding the optimum for the local copies of the variables, subject to the condition that all local copies are equal to the global variable. The problem is solved iteratively, with all local copies eventually converging to the global optimal value.

For the problem described by Eqs.~\eqref{LinDistFlow}, each node $j$ will keep the local copies of the $Q_j$, $Q_{j-1}$, $U_j$ and $U_{j-1}$ variables, which we will denote $Q^{+}_j$, $Q^{-}_j$, $U^{+}_{j}$ and $U^{-}_{j}$  respectively. We note that the algorithm assumes that one can measure power flows $P_j$ between nodes which are treated as auxiliary constant parameters of the algorithm.
The optimization problem \eqref{LinDistFlow} formulated as a consensus problem becomes:
\begin{subequations}\label{LinDistFlowconsensus}
\begin{eqnarray}
%\label{LinDistFlowQconsensus}
&\min\limits_{Q}& \sum_{j=1}^{n} \frac{r_{j-1} \left(Q^{-}_j \right)^2 }{V_0^2}, \label{LossLinearQconsensus} \\
&\mbox{s.t.}&  \forall j=1,\ldots,n: \nonumber \\
&& Q^{+}_j - Q^{-}_{j} + q_j^{(c)}  - \tilde{s}_j \leq 0, \label{PV_constraint_1con} \\
&&  -Q^{+}_j + Q^{-}_{j}  - q_j^{(c)}  - \tilde{s}_j \leq 0,  \label{PV_constraint_2con} \\
&& Q^{+}_j = Q_j, \quad  Q^{-}_j = Q_{j-1}, \quad Q_n = 0  \label{consensus} \\
&&  V_0^2 \left(\epsilon^2 - 2\epsilon \right) \leq U^{+}_j\leq V_0^2\left(\epsilon^2 + 2\epsilon\right),\label{Voltage-Uconsensus} \\
%&&\forall j = 0,\ldots,n-1: \nonumber \\
&& U^{+}_{j} = U^{-}_j - 2(r_{j-1} P_{j-1}+x_{j-1} Q^{-}_{j}), \label{V-linconsensus} \\
&&  U^{+}_{j} =  U_j,\quad U^{-}_j = U_{j-1} ,\quad U^{-}_1 = 0 \label{V-consensus}
\end{eqnarray}
\end{subequations}
The conditions \mref{consensus,V-consensus} ensure that all local copies of the variables are equal to the global variables $Q_j$ and $U_j$
and hence the optimization problem \eqref{LinDistFlowconsensus} is equivalent to \eqref{LinDistFlow}.
We solve \eqref{LinDistFlowconsensus} using distributed ADMM method \cite{Boyd2011,boyd2013}, for which
the augmented Lagrangian is
\begin{subequations}
% \begin{eqnarray}
%   {\cal L}^{\rm{ADMM}} \left(Q^{+},Q^{-},Q,U^{+},U^{-},U,\lambda^{Q^{+}}, \lambda^{Q^{-}}, \lambda^{U^{+}}, \lambda^{U^{-}}, \rho \right) = & \\  \sum_{j=1}^{n} {\cal L}^{\rm{ADMM}}_j, & \nonumber
%
% \end{eqnarray}
\begin{equation}
  {\cal L}^{\rm{ADMM}} =  \sum_{j=1}^{n} {\cal L}^{\rm{ADMM}}_j,
  \end{equation}
where
\begin{eqnarray}
 {\cal L}^{\rm{ADMM}}_j  &=& \frac{r_{j-1} \left( Q^{{-}}_j \right)^2} {V_0^2}  \label{LagrangianADMM} \\
 &+& \frac{\rho}{2} \left( Q^{+}_j - Q_j \right)^2 + \frac{\rho}{2} \left( Q^{ -}_{j} - Q_{j-1}  \right)^2 \nonumber \\
  &+& \frac{\rho}{2} \left( U^{+}_j - U_j \right)^2 + \frac{\rho}{2} \left( U^{-}_{j} - U_{j-1}  \right)^2  \nonumber \\
 &+& \lambda^{Q^+}_j \left( Q^{+}_j - Q_j  \right) +  \lambda^{Q^{-}}_j \left( Q^{-}_j - Q_{j-1}  \right), \nonumber \\
  &+& \lambda^{U^+}_j \left( U^{+}_j - U_j  \right) +  \lambda^{U^{-}}_j \left( U^{-}_j - U_{j-1}  \right). \nonumber
\end{eqnarray}
\end{subequations}
The quadratic terms in the objective function with $\rho/2$ prefactor represent penalties for the local variables being different from the global variables. These terms do not change the optimal value, as the constraints \mref{consensus,V-consensus} require that the local and  global variables are equal at the optimum. The dual variables associated with \mref{consensus,V-consensus} are $\lambda^{Q^{+}}$, $\lambda^{Q^{-}}$, $\lambda^{U^+}_j$ and $\lambda^{U^{-}}_j$. Note that we do not include constraints \mref{PV_constraint_1con,PV_constraint_2con} in the Lagrangian, as the algorithm will be minimizing ${\cal L}_{\rm{ADMM}}$ in such a way that $Q^{+}$ and $Q^{-}$ will always stay within the feasible set, i.e. satisfy \mref{PV_constraint_1con,PV_constraint_2con}.

The ADMM distributed consensus algorithm is an iterative algorithm where the $k+1$-th iteration starts with values $Q_j(k)$, $Q_{j-1}(k)$, $Q_j^{+}(k)$,$Q_j^{-}(k)$, $\lambda^{Q^{+}}_j(k)$, $\lambda^{Q^{-}}_j(k)$, $U_j(k)$, $U_{j-1}(k)$, $U_j^{+}(k)$,$U_j^{-}(k)$, $\lambda^{U^{+}}_j(k)$, $\lambda^{U^{-}}_j(k)$ for each node $j$. One iteration of the algorithm consists of the following steps:
\begin{enumerate}
 \item  {\bf Minimization step}. For each node $j$, the following optimization problem is solved
 \begin{eqnarray}
   && \min\limits_{Q^{-}_j,Q^{+}_j, U^{+}_j, U^{-}_j}  {\cal L}^{\rm{ADMM}}_j   \label{LocalLoss} \\
   && \mbox{s.t. Eqs.~\mref{PV_constraint_1con,PV_constraint_2con,Voltage-Uconsensus,V-linconsensus}} \nonumber
 \end{eqnarray}
This minimization step is a convex optimization problem with quadratic objective function of four local variables $Q^{-}_j,Q^{+}_j, U^{-}_j, U^{+}_j$ with linear constraints and can be solved analytically by evaluating the corresponding Karush-Kuhn-Tucker conditions \cite{boyd2004convex}. However, the expressions are bulky and we do not present them here for the sake of brevity.
 Each node $j$ can perform the minimization step independently, as the optimization is carried over the local variables at node $j$.  The solutions to the local minimization problem in $k$-th iteration of the ADMM distributed consensus algorithm are denoted as $Q^{-}_j(k+1)$, $Q^{+}_j(k+1)$,
  $U^{-}_j(k+1)$, $U^{+}_j(k+1)$.
 \item {\bf Averaging step}. This step updates the global (shared) variables $Q$ and $U$. The update rules for each are the following:
 \begin{eqnarray}
  && \forall j = 1,\ldots,n-1: \nonumber \\
  && Q_j(k+1) = \frac{1}{2} \left( Q^{+}_j(k+1) + Q^{-}_{j+1}(k+1)\right) \nonumber \\
  && U_j(k+1) = \frac{1}{2} \left( U^{+}_j(k+1) + U^{-}_{j+1}(k+1)\right) \nonumber \\
   && Q_n(k+1) = 0, \, Q_0(k+1) =  Q^{-}_1(k+1),\, U_n = U^{+}_n  \nonumber
 \end{eqnarray}
 This step requires communication between nearest neighbors, as they need to exchange their local variables in order for each node $j$ to calculate the new value for $Q_j$ and $Q_{j-1}$ which is the average of the respective local copies of neighboring nodes.

 \item {\bf Lagrange multipliers update step}.
 The Lagrange multipliers, which are also stored by each node locally, are updated according to the following rules for each node $j$:
 \begin{eqnarray*}
  {\lambda^{Q^{+}}_j}(k+1) &=&  {\lambda^{Q^{+}}_j}(k) + \rho \left(  Q^{+}_j(k+1) - Q_j(k+1) \right) \\
  {\lambda^{Q^{-}}_j}(k+1) &=&  {\lambda^{Q^{-}}_j}(k) + \rho \left(  Q^{-}_j(k+1) - Q_{j-1}(k+1) \right)\\
   {\lambda^{U^+}_j}(k+1) &=&  {\lambda^{U^+}_j}(k) + \rho \left(  U^{+}_j(k+1) - U_j(k+1) \right) \nonumber\\
  {\lambda^{U^{-}}_j}(k+1) &=&  {\lambda^{U^-}_j}(k) + \rho \left(  U^{-}_j(k+1) - U_{j-1}(k+1) \right) \nonumber.
 \end{eqnarray*}
 All variables involved in this step have been calculated and communicated in the previous step, which means that the Lagrange multipliers are updated locally at each node $j$.
\end{enumerate}

%Each step of the ADMM consensus algorithm is more intense computationally than an individual step of the dual ascent algorithm: ADMM
%Every step of the ADMM consensus algorithm
%requires each of the nodes to solve a constrained quadratic programming problem with four variables.

The ADMM algorithm requires synchronized communication between the neighboring nodes where local variables ($Q^{+}_j$, $Q^{-}_j$, $U^{+}_j$ and $U^{-}_j$) are communicated between nearest neighbors.
These local variables can be interpreted as 'beliefs' of node $j$ about which reactive power should be flowing in and out of the node and what should be the voltage magnitudes. The consensus algorithm leads to convergence of these local variables between the neighboring nodes, thus finding a global optimal solution.
Once the algorithm converges, the local variables will actually correspond to an optimized feasible solution $Q_j$ and $U_j$ of problem \eqref{LinDistFlowconsensus}. The actual values of reactive power injected by inverters can be calculated by each node from its local variables as
\begin{equation}
 q^{(g)}_j =  Q^{+}_j - Q^{-}_j  + q^{(c)}_j
 \label{how_to_q_g}
\end{equation}
The above solution for $q^{(g)}_j$ is guaranteed to be within the allowed bounds given by the $\tilde{s}_j$, as the local variables always satisfy the conditions \mref{PV_constraint_1con,PV_constraint_2con}.

%Hence even before the algorithm converges to the optimal solution, the sub-optimal solution for $q_j^{(g)}$ can be obtained from \eqref{how_to_q_g}.
%as opposed to the dual-ascent algorithm, where prior to the variables $Q$ convergence, they do not actually have to satisfy the conditions \eqref{PV_constraint_1A} and \eqref{PV_constraint_1B}.

To test the performance of the ADMM algorithm, we compare its convergence with a dual ascent algorithm, which we derive for a simplified LinDistFlow problem without voltage constraints in Section \ref{sec_noU}.
As we will show in Section \ref{sec_cases}, the ADMM distributed consensus algorithm converges faster than the dual ascent algorithm, which is known to be its general main advantage \cite{Boyd2011}.

\subsection{Dual ascent algorithm for distributed control of the inverters with no voltage constraints}
\label{sec_noU}
We now formulate a dual ascent approach to solving the simplified optimization problem \eqref{LinDistFlow} without considering the voltage constraints \mref{V-lin,Voltage-U}.
%The objective function in \eqref{LinDistFlowQnoU} is separable and the constraints \eqref{PV_constraint_1A} and \eqref{PV_constraint_1B} are local.  The conditions admit a message passing approach to the solution.
%We first discuss the dual ascent algorithm \cite{bertsekas1999nonlinear}.
The Lagrangian becomes
\begin{subequations}
\begin{eqnarray}
{\cal L}\left(Q,\zeta^+,\zeta^- \right) &=& \sum_{j=0}^{n-1}\Biggl[\frac{r_j Q_j^2}{V_0^2}  +  \\
  &+& \zeta^+_j\left(Q_{j+1} - Q_j -\tilde{s}_{j+1}+q_{j+1}^{(c)}\right)  \nonumber \\
  &+& \zeta^-_j\left(-Q_{j+1}+Q_j-\tilde{s}_{j+1}-q_{j+1}^{(c)}\right)  \nonumber  \Biggr].
\label{Lagrangian}
\end{eqnarray}
\end{subequations}
The dual ascent algorithm consists of the following steps:
\begin{subequations}
\begin{enumerate}
 \item Minimize ${\cal L}\left(Q,\zeta^+,\zeta^- \right)$ over $Q$ for given $\zeta^+, \zeta^-$, which
% \begin{equation}
%   \forall j = 0,\ldots,n-1 : \quad \frac{\partial {\cal L}}{\partial Q_j} = 0.
% \end{equation}
leads to the following update rule in $k$-th iteration
 \begin{equation}
 \label{ascent_minimization}
  Q_j (k+1) = \frac{V_0^2}{2 r_j} \left( \zeta^{+}_j(k) - \zeta^{+}_{j-1}(k) + \zeta^{-}_{j-1}(k) - \zeta^{-}_{j}(k) \right).
 \end{equation}

 \item Update dual variables according to
 \begin{eqnarray}
  \zeta^{+}_j(k+1) &=&  \max  \left(0, \zeta^{+}_j(k) + \alpha  \Delta^+_{j+1} \right) \label{ascent_alphaplus} \\
  \zeta^{-}_j(k+1) &=&  \max \left(0, \zeta^{-}_j(k) + \alpha  \Delta^-_{j+1} \right)  \label{ascent_alphaminus}
 \end{eqnarray}
\end{enumerate}
where
\begin{eqnarray*}
 \Delta^+_{j+1} &=&   Q_{j+1}(k+1) - Q_j(k+1) + q_{j+1}^{(c)} - \tilde{s}_{j+1}  \\
 \Delta^-_{j+1}&=&  - Q_{j+1}(k+1) + Q_j(k+1) - q_{j+1}^{(c)} - \tilde{s}_{j+1}.
\end{eqnarray*}
\end{subequations}
This scheme allows parallel implementation, where each node $j$ receives values of $\zeta^{+}_{j}(k)$ , $\zeta^{-}_{j}(k)$ from its right neighbor $j+1$ and sends values of
$\zeta^{+}_{j-1}(k)$ , $\zeta^{-}_{j-1}(k)$ to its left neighbor $j-1$. Node $j$ then calculates $Q_{j}(k+1)$ using those variables and then sends the result to neighbor $j+1$, while receiving $Q_{j-1}(k+1)$ from neighbor $j-1$. The communicated values are then used by node $j$ to calculate $\zeta^{+}_{j-1}(k+1)$ and $\zeta^{-}_{j-1}(k+1)$.

An advantage of the dual ascent algorithm is that it requires the nodes to perform trivial algebraic operations (which are simpler than the solution of Eq.~\eqref{LocalLoss}) and synchronously communicate their local variables $Q$, $\zeta^+$, and $\zeta^-$ to their neighbors.
However, based on discussion in \cite{Boyd2011}, the dual ascent method is expected to require a large number of iterations to converge.  The speed of convergence is controlled by parameter $\alpha$, but the range of feasible $\alpha$ is limited---choosing $\alpha$ which is too large results in a failure to converge, while $\alpha$ chosen too small translates into a slow convergence.

We provide comparison of the dual ascent algorithm with the ADMM algorithm in Section \ref{sec_results}. Since the proposed dual ascent algorithm does not
consider voltage constraints, we will also consider a version of the ADMM algorithm without the voltage constraints (referred to as ADMM-noV), which can be straightforwardly obtained from the ADMM algorithm by excluding voltage variables $U_j, U^+_j, U_j^-$ and removing constraints \mref{Voltage-Uconsensus,V-linconsensus,V-consensus}.

\section{Experiments with distributed global optimization on different distribution circuit cases}
\label{sec_cases}

We explore the performance of our distributed optimization algorithms via simulations on a range of distribution circuit cases. We first introduce the considered feeder line configuration cases and then compare the performance of the dual ascent, ADMM and ADMM-noV algorithms on these cases. The global optimization
results from these algorithms are compared with the sub-optimal local optimization scheme proposed in \cite{KostyaMishaPetrScott}, where each node only uses its local information about $q_j^{(c)}$  to set its $q_j^{(g)}$. Finally, we compare losses and voltages calculated with the LinDistFlow equations with the ones calculated by DistFlow equations for the same set of injected reactive power $q^{(g)}$ in order to check the validity of the underlying approximation.

\subsection{Distribution circuit test cases}

The properties of the distribution circuit used in the simulations are summarized in the following table:
\begin{center}
 \begin{tabular}{|c|cclll|}
  \hline
  {\bf Case} & {\bf Nodes } &  {\bf PV-pen} & $\mathbf{p^{(c)}_{\rm{\bf max}}}$ & $\mathbf{p^{(g)}}$ & $\mathbf{s_{\rm{\bf max}}}$ \\
  \hline
    1 & 100 & 100\% &  $4$ kW   & $1$ kW & $1.1$ kW \\
    2 & 100 & 50\%  &  $4$ kW   & $1$ kW & $1.1$ kW \\
    3 & 250 & 50\%  &  $2.5$ kW & $1$ kW & $2.2$ kW \\
    4 & 250 & 50\%  &  $1$ kW   & $2$ kW & $2.2$ kW \\
    5 & 150 & 85\%  &  $4$ kW &  $0.9$ kW & $1.1$  kW \\
    6 & 200 & 100\% &  $3.75$ kW & $0$ kW & $2.2$ kW \\
    %7 & 150 & 70\%  &  $2$ kW &  $6.5$ kW & $10$  kW \\  \hline
    7 & 150 & 70\%  &  $2$ kW &  $7.0$ kW & $10$  kW \\  \hline

 \end{tabular}
\end{center}

The ``PV-pen'' column indicates the percentage of nodes in the distribution circuit that have PV generation installed.  These nodes all inject the same power which is given in the column denoted as $p^{(g)}$. $s_{\rm{max}}$ is the apparent power capacity of the inverters which enters into the constraints in Eq.~\eqref{PV_constraint}. We set $s_j = s_{\rm{max}}$ for all the nodes with PV generation installed, and $s_j = 0$ for the rest.  The real power consumed at each node, $p_j^{(c)}$, is chosen from a uniform distribution between $0$ and $p^{(c)}_{\rm{max}}$. The reactive power consumption is set to $q^{(c)}_j = 0.25  p_j^{(c)}$ for all the cases except for cases 5 and 7 where $q^{(c)}_j = f_j p_j^{(c)}$ with $f_j$ is drawn from a uniform distribution between $0.01$ and $1.0$ for case 5 and between $0$ and $1.0$ for case 7. The distribution circuit line resistance $r_j$ and reactance $x_j$ are $0.33\, \Omega/\rm{km}$ and $0.38\, \Omega / \rm{km}$ respectively, with the distances between the neighbors always set to
$0.25\,\rm{km}$. The voltage at the start of the distribution circuit is $V_0 = 7.2\,{\rm kV}$.

The considered distribution circuits are based on feeder line configurations used in our previous work \cite{KostyaMishaPetrScott,KostyaMishaPetrScott2, KostyaMishaPetrScott3} where suboptimal, policy-based control schemes  were analyzed. Case 6 corresponds to higher loads with no generation but 100\% penetration (representing example of a nighttime case with high loads).  This case is included because the node voltages violate the {\it minimum} voltage constraint without control of the $q^{(g)}$ (i.e., when all $q_j^{(g)} = 0$) . Case 7 corresponds to a scheme with high power generation and low consumption in the feeder line such that, without control, the voltage violates the {\it maximum} voltage constraint. The voltage rise above the allowed limit is a possible issue in networks with high distributed PV generation   
\cite{tonkoski2012impact}.

\subsection{Global vs local optimization}
\label{sec:global_vs_local}
First, we compare the globally optimal solutions for the injected reactive power obtained with the distributed algorithms described above (denoted as $\tilde{q}_j^{(g)}$) with the local policy-based scheme proposed in \cite{KostyaMishaPetrScott}. In the local scheme, the inverters are set to supply the local reactive power consumption up to their apparent power capacity:
\begin{equation}
\label{q_local}
 q^{(g)}_j({\rm{local}}) = \min \left( q^{(c)}_j , \tilde{s}_j \right).
\end{equation}

We compute the total circuit losses (using Eq.~\eqref{LossLinear}) for both global and local policy-based scheme (setting reactive powers to  $\tilde{q}_j^{(g)}$  and $q^{(g)}_j({\rm{local}})$ respectively) and divide by the losses for the ``no-optimization scheme,'' i.e. with all $q^{(g)}_j$ set to zero. This scheme corresponds to the current situation for PV inverters which inject real power at a nominal power factor of $1.0$.  The normalized losses from the globally-optimal solution of Eqs.~\eqref{LinDistFlow} (Loss$^{\rm glob}_{\rm lin}$) and the local policy control (Loss$^{\rm loc}_{\rm lin}$) are shown in the following table for cases 1-5: \\
\begin{center}
\begin{tabular}{|c|cc|}
\hline
 {\bf Case} & {\bf Loss$^{\rm glob}_{\rm lin}$} & {\bf {\bf Loss$^{\rm loc}_{\rm lin}$}} \\
 \hline 1 & $0.834$ & $ 0.845$ \\
 2 &  $0.941$ & $0.949$ \\
 3 &  $0.847$ & $0.890$ \\
 4 &  $0.954$ &  $0.962$ \\
% 5 &  $0.941255$ & $0.941176$ \\
 5 &    $0.700$  &    $0.771$           \\  \hline
\end{tabular}
\end{center}
In each case, both control schemes lower the total loss of real power, and except for cases 3 and 5, the local policy-based scheme performs nearly as well as the globally optimal solution.  This should be expected because these cases are not so heavily loaded or over-generated that the voltage exceeds the normal operational limits.  As was shown in \cite{kundu2013distributed}, the local policy in \eqref{q_local} is approximately an optimal policy until these limits are approached.  However, this local policy is unable to respond appropriately when the voltage deviates beyond its normal operational limits.

Cases 6 and 7 are shown in the table below, along with the minimum (for case 6) and maximum (for case 7) voltages normalized with respect to its nominal value $V_0$.
The voltages shown in the table were calculated by solving the exact DistFlow eqs.~\eqref{DistFlow} using the values of $q_j^{(g)} = \tilde{q}_j^{(g)}$ for the global optimization of LinDistFlow eqs.~\eqref{LinDistFlow} scheme (V$^{\rm glob}$), $q_j^{(g)} = q^{(g)}_j({\rm{local}})$ for local policy control scheme (V$^{\rm loc}$), and $q_j^{(g)} = 0$ for the no-optimization scheme (V$^{\rm{noopt}}$) respectively.

The voltages in the no-optimization scheme violate the nominal operating limits for both cases. The local policy-based control is unable to fully correct this situation and, in the case of voltage in case 7, makes the situation worse.  The global optimization scheme enables the voltage to be corrected to respect the constraints, which naturally leads to higher losses. We note that the minimum (maximum) normalized voltage calculated with LinDistFlow approximation is $0.95$ ($1.05$) for case 6 (case 7) for global optimization scheme, but the actual values obtained from DisFlow eqs.~\eqref{DistFlow} are slightly below the voltages calculated from the linear approximation \eqref{LinDistFlow}. The agreement between the linear approximation and the DistFlow equations will be further discussed in Section \ref{sec:validity}.
%(ADMM-V is required in this case to correct the voltage deviations)
\begin{center}
\begin{tabular}{|c|ccccc|}
 \hline {\bf Case} & {\bf Loss$^{\rm glob}_{\rm lin}$} & {\bf Loss$^{\rm loc}_{\rm lin}$}& {\bf V}$^{\rm glob}$ & {\bf V}$^{\rm loc}$ & {\bf V}$^{\rm{noopt}}$ \\ \hline
   % 6 &   $0.954$ &  $0.941$  & $ 0.950 $  &  $0.941$  & $ 0.923$ \\  %calculated with lin
    6 &   $0.954$ &  $0.941$  & $ 0.947 $  &  $0.938$  & $ 0.920$ \\  %calculated with lin
    7 &   $1.11$ & $1.09$  & $1.045$ & $1.074$ & $1.071$ \\ \hline  %calculated with nonlin
\end{tabular}
\end{center}

\subsection{Performance of distributed optimization algorithms}
\label{sec_results}

\begin{figure}
\centering
\includegraphics[width=8.5cm]{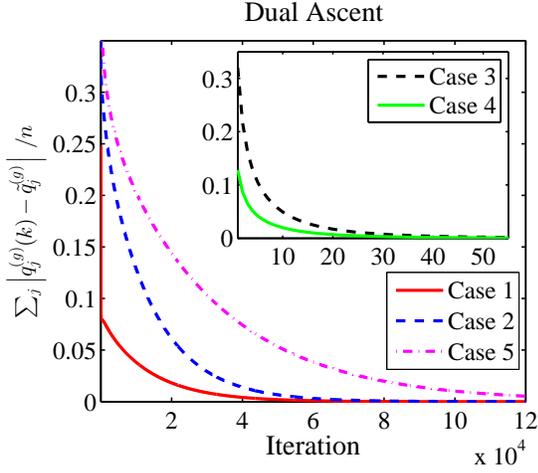}\\~\\
\caption{
The average absolute deviation from the optimal values for $q^{(g)}$  as a function of number of iterations for the dual ascent algorithm applied to cases 1-5. The cases 3 and 4 are shown in the inset, while the cases 1, 2 and 5 are plotted with the iteration axis scaled by the multiples of $10^4$ iterations.
}
\label{fig_ascent}
\end{figure}

\begin{figure}
\centering
\includegraphics[width=8.5cm]{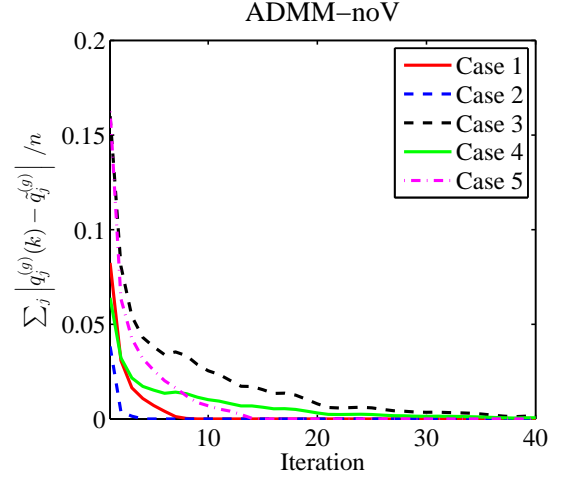}\\~\\
\caption{
The average absolute deviation from the optimal values for $q^{(g)}$ as a function of number of iterations. The data were obtained from ADMM-noV algorithm for cases 1-5.
}
\label{fig_admm}
\end{figure}

\begin{figure}
\centering
\includegraphics[width=8.5cm]{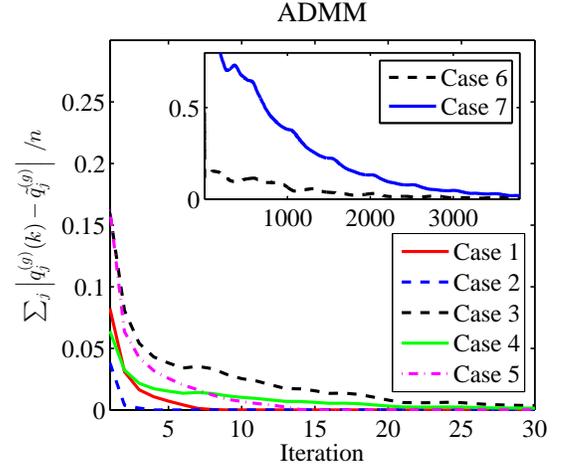}\\~\\
\caption{
The average absolute deviation from the optimal values for $q^{(g)}$  as a function of number of iterations. The data were obtained from ADMM algorithm for cases 1-7, with cases 6 and 7 shown separately in the inset.
}
\label{fig_admm_v}
\end{figure}

For cases 1-5, where voltage constraints are satisfied even without explicitly including them in the optimization problem, the dual ascent, ADMM-noV and ADMM algorithms yield identical minimum real power losses.  However, these algorithms do not have equal performance in terms of convergence.
For  cases 1-5, we contrast the convergence properties of each algorithm by comparing the intermediate values of $q^{(g)}_j$ at iteration $k$ with the global optimum $\tilde{q}^{(g)}_j$ obtained by solving \eqref{LinDistFlow} with the CVX package \cite{cvx,gb08}. We plot on the y-axis the average absolute deviation as a function of number of iterations, defined as
\begin{equation}
 {\rm{D}}(k) = \frac{1}{n} \sum_{j=1,\ldots,n} \left| q^{(g)}_j(k) - \tilde{q}^{(g)}_j \right|
\end{equation}
where $n$ is the total number of nodes in the feeder line.  At each iteration, the values of $q^{(g)}(k)$ are calculated from the local variables. For the dual ascent algorithm, we utilize $Q(k)$ and \eqref{Q-lin} to obtain $q^{(g)}(k)$ while for ADMM-noV and ADMM we use the variables $Q^{-}(k)$, $Q^{+}(k)$ and Eq.~\eqref{how_to_q_g}.

The values of ${\rm{D}}(k)$ for dual ascent, ADMM-noV, and ADMM algorithms are plotted in Figs.~\ref{fig_ascent}, \ref{fig_admm} and \ref{fig_admm_v}, respectively. For all the algorithms, the initial value of $Q_j$ and of all Lagrange multipliers was chosen to be equal to $0$. The initial values of $U_j$ for the ADMM algorithm were taken from the solution of Eq.~\eqref{V-lin}, with $Q_j$ set to 0. %The other choices of initial values for ADMM are discussed in appendix \ref{sec_initial_point}.

For the cases 1, 2, and 5, the dual ascent algorithm (Fig.~\ref{fig_ascent}) takes the order of $10^4$ iterations to converge while cases 3 and 4 only required a few tens of iterations.  The fast convergence is primarily due to the optimal solution for those cases being close to the initial choice of $Q_j = 0$. The speed of convergence is controlled by the parameter $\alpha$ in Eqs.~\mref{ascent_alphaplus,ascent_alphaminus}. Via empirical experimentation, we find $\alpha = 0.05/V_0^2$ to give optimal convergence performance as choosing larger/smaller $\alpha$ caused numerical instability/slower convergence.

For cases 1-5, the ADMM-noV (Fig.~\ref{fig_admm}) and ADMM (Fig.~\ref{fig_admm_v}) algorithms converge to the optimal $q^{(g)}$ within tens of iterations. This rapid convergence dramatically outperforms the dual ascent algorithm in cases where the initial guess is not close to the optimal solution.  In the cases 6 and 7, ADMM requires order of $10^3$ iterations (Fig.~\ref{fig_admm_v}, inset).  The reduced performance in this case is because the voltage constraints are violated at some nodes, and this information needs to propagate throughout the entire distribution circuit.

In general, the dual ascent algorithm requires more iterations, and hence more rounds of communication between nodes, to converge. The adjustment of the $q^{(g)}$ should be carried out on the faster timescales than the $p^{(g)}$ and $p^{(c)}$ are changing.  For PV generation, these changes can be on the order of one to several minutes.  In some cases, the large number of interactions required for dual ascent would challenge the capability of grid communications systems.

\begin{figure}
\centering
\includegraphics[width=8.5cm]{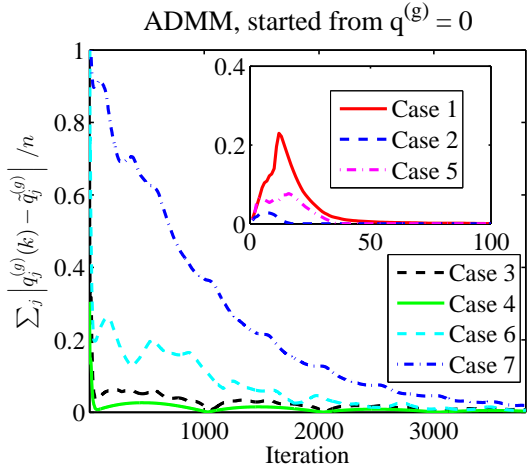}\\~\\
\caption{
The average absolute deviation from the optimal values for $q^{(g)}$  as a function of number of iterations. The data were obtained from ADMM algorithm for cases 1-7.
The initial state was taken for $Q,\,Q^{+},\,Q^{-}$ and $U$ corresponding to solution of power flow equations when $q^{(g)}_j = 0$.
}
\label{fig_admmv_fromqg0}
\end{figure}

\begin{figure}
\centering
\includegraphics[width=8.5cm]{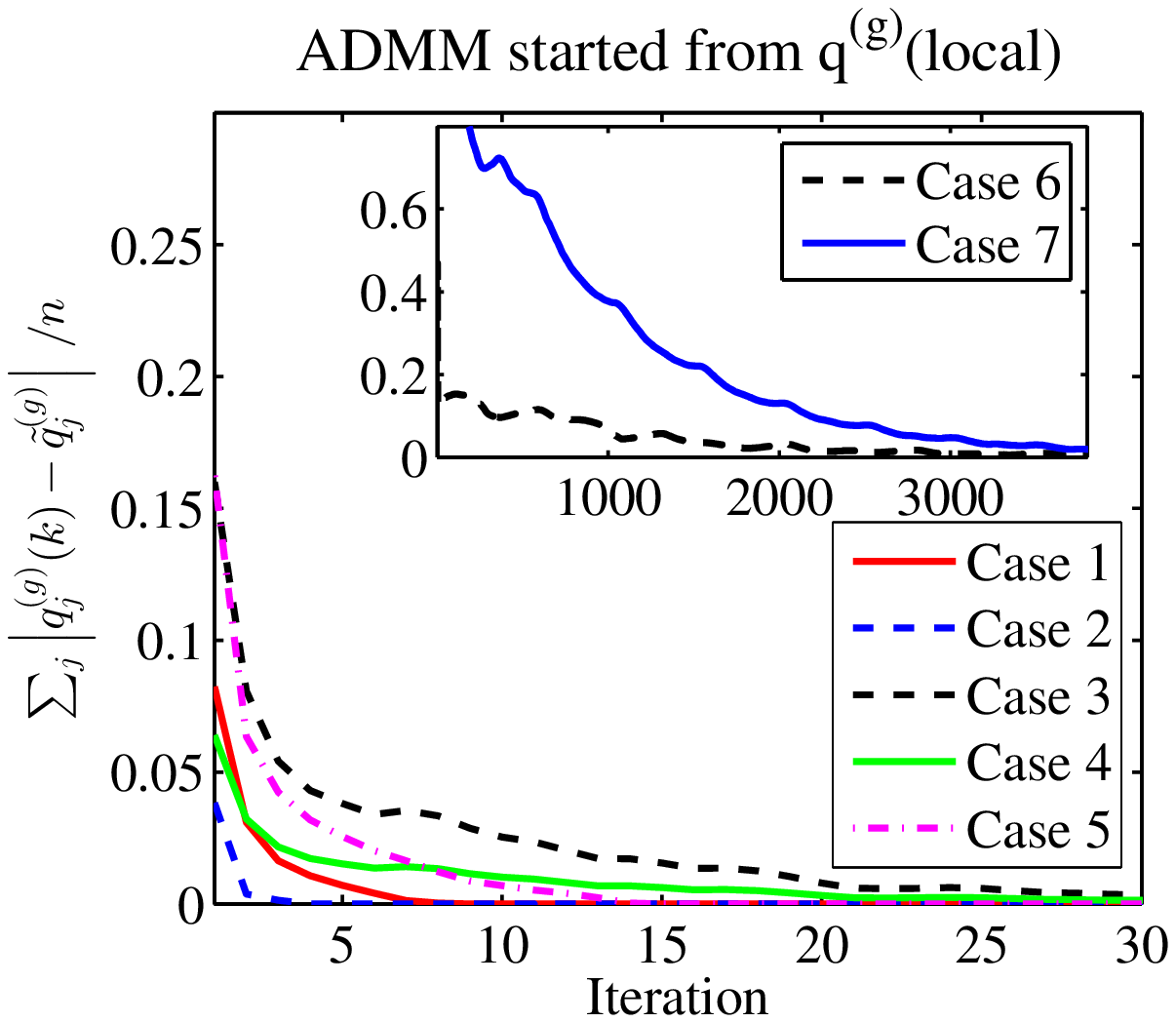}\\~\\
\caption{
The average absolute deviation from the optimal values for $q^{(g)}$  as a function of number of iterations. The data were obtained from ADMM algorithm for cases 1-7.
The initial state was taken for $Q,\,Q^{+},\,Q^{-}$ and $U$ corresponding to solution of power flow equations when $q^{(g)}_j = q^{(g)}_j(\rm{local})$.
}
\label{fig_admmv_fromqlocal}
\end{figure}

\begin{figure}
\centering
\includegraphics[width=8.5cm]{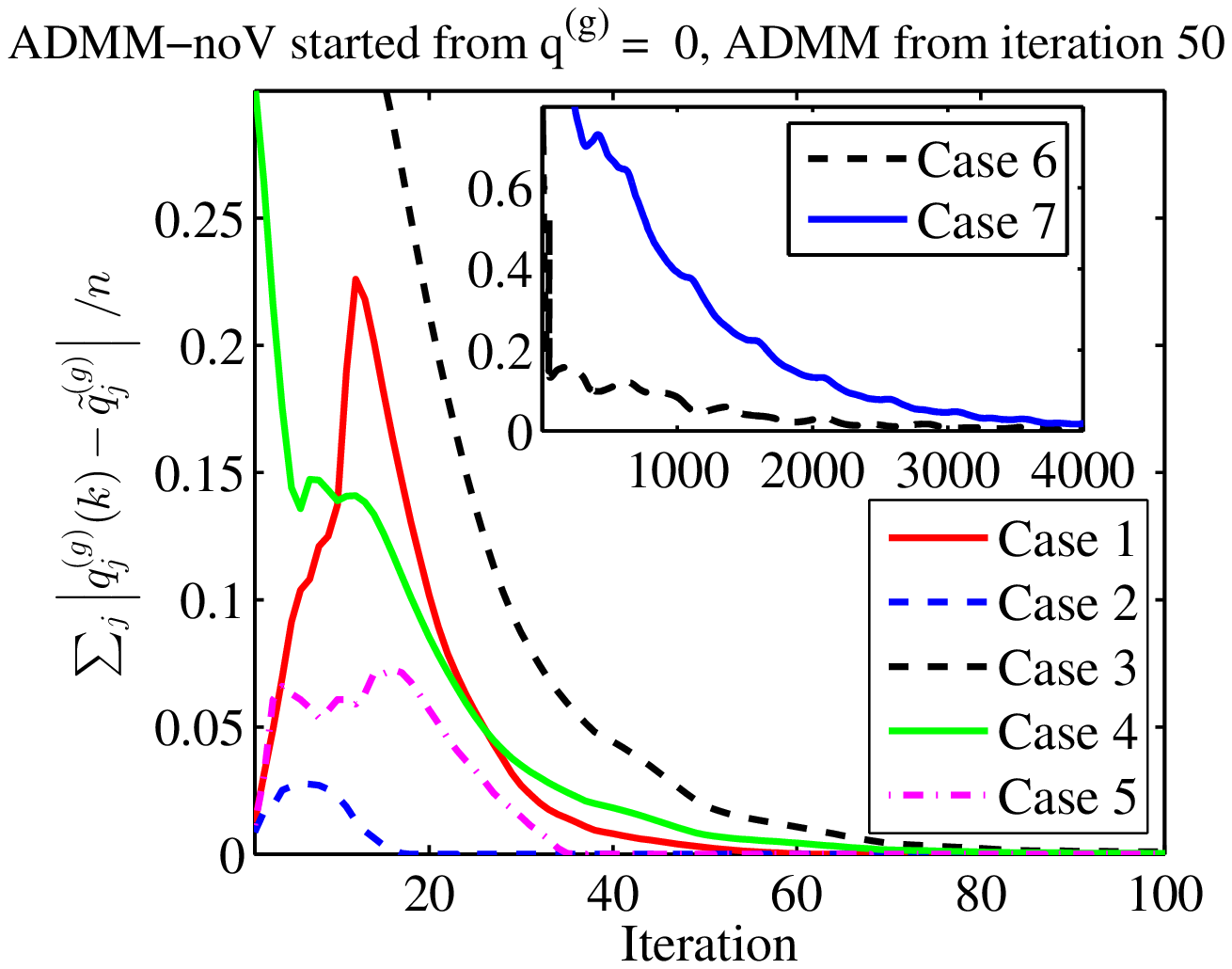}\\~\\
\caption{
The average absolute deviation from the optimal values for $q^{(g)}$  as a function of number of iterations. The data were obtained from an algorithm which first ran ADMM-noV for 50 iterations and then used the measured values of voltage to initialize ADMM algorithm.
The initial state was taken for $Q,\,Q^{+},\,Q^{-}$ corresponding to solution of power flow equations when $q^{(g)}_j = 0$.
}
\label{fig_admmv_fromqlocalUafter}
\end{figure}

Our choice of initial condition in for Figs.~\ref{fig_admm}, \ref{fig_admm_v} was simple $Q = Q^+=Q^-=0$, $U$ was set in accordance with solutions of \eqref{V-lin} with $Q = 0$.
%, and, finally, all the Lagrangian multipliers were also set to zero.
Such choice of $Q$ allowed for direct comparison between the dual ascent and ADMM algorithms, as they started from the same point. Note that the choice of nonzero initial  $Q_j$ in the dual ascent algorithm would also require setting up nonzero Lagrange multipliers $\zeta^{+}_j, \zeta^{-}_j$, according to \eqref{ascent_minimization}, that would likely require additional (preemptive) communication between the nodes to agree on the multipliers.

The zero initial choice of $Q$ was sufficiently close to the actual optimal solutions in cases 3 and 4. However, other cases required more iterations to converge. On the other hand, the ADMM-noV algorithm converged within tens of iterations even in the other cases (1, 2 and 5) in spite of the fact that the final solutions for $Q$ were rather different from the zero initial guess. The ADMM algorithm, applied to the cases 6 and 7 where voltage constraints were violated in the unoptimized solution, took much longer ($O(10^3)$) to converge. Experimenting with the ADMM we have found that its convergence is rather sensitive to the initial guess of voltages. In particular, we tested the following two initializations:
(a) $Q$ and $U$ corresponding to the case where $q^{(g)}_j = 0$ for all nodes (the unoptimized solution), and (b) $Q$ and $U$ corresponding to the state where  $q^{(g)}_j = q^{(g)}_j(\rm{local})$, where $q^{(g)}_j(\rm{local})$ is the local optimum defined by Eq.~\eqref{q_local}. Convergence of the ADMM is compared for these two initial states in Figs.~\ref{fig_admmv_fromqg0} and \ref{fig_admmv_fromqlocal} respectively.  Note that the bare ADMM-noV (i.e. no voltage constraints) converged in tens of iterations for both (a) and (b) initiations. Convergence of ADMM initialized with $q^{(g)}_j=q^{(g)}_j(\rm{local})$ (Fig.~\ref{fig_admmv_fromqlocal}) matched convergence of the bare ADMM-noV for cases 1-5, as does ADMM when initialized from $Q_j = 0$ (Fig.~\ref{fig_admm_v}).  However, the same ADMM initialized with $q^{(g)}_j = 0$  took $O(10^3)$ iterations to converge in the cases 3 and 4,  presumably due to rather inaccurate voltage guess/initiation. One possible way to get around this problem is to run the ADMM-noV algorithm for few iterations, and then switch to the ADMM initializing voltages from the measured values at the nodes, as shown in Fig.~\ref{fig_admmv_fromqlocalUafter}.

\begin{figure}
\centering
\includegraphics[width=8.5cm]{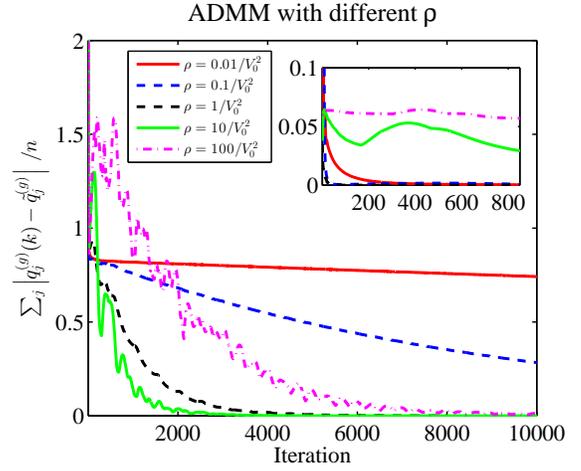}\\~\\
\caption{
The average absolute deviation from the optimal values for $q^{(g)}$ for case 7 (in the main plot) and case 4 (in the inset) for different values of parameter $\rho$.
The initial state was taken for $Q,\,Q^{+},\,Q^{-}$ and $U$ corresponding to solution of power flow equations when $q^{(g)}_j = q^{(g)}_j(\rm{local})$.
}
\label{fig_admmv_rhos}
\end{figure}

The convergence of ADMM algorithms is further influenced by the choice of the parameter $\rho$ in the ADMM Lagrangian \eqref{LagrangianADMM}. We used $\rho = 1/V_0^2$, which we found, via empirical experimentation, to give the best convergence properties for our algorithm.
We have considered decreasing and increasing $\rho$ by a factor of 10 and 100. The choice which provides good convergence properties for all the cases and initial conditions considered is found to be  $\rho = 1/V_0^2$, even though for some particular cases and choice of initial conditions, smaller or larger $\rho$ sometimes performs slightly better. The convergence rate is sensitive to the choice of $\rho$, as choosing too large or too small $\rho$ can lead to significant slowdown of convergence in some cases, as illustrated in Fig.~\ref{fig_admmv_rhos}.

Finally, we note that for the numerical implementation, the actual units of representation of voltage and reactive power can play significant role, as they balance the magnitude of the terms in the ADMM Lagrangian.
For our ADMM simulations, we represented $Q$ in kW and $U$ in $100$ kV$^2$ units. Choosing smaller magnitude for the representation of the voltages could lead to accumulation of errors due to the differences between the local variables $U^{+}, U^{-}$ and the actual value of $U$, which would then affect the convergence for algorithms where voltage constraints are violated typically at the last few nodes in the distribution circuit. Alternatively, one could consider having two different parameters $\rho_Q$ and $\rho_U$ and adapt them separately for optimal convergence.

\subsection{Validity of LinDistFlow approximation}
\label{sec:validity}
The optimization algorithms formulated above rely on the validity of LinDistFlow equations. To test the accuracy of this approximation, we take the $\tilde{q}_j^{(g)}$ solution of the global optimization formulated with the LinDistFlow equations \eqref{LinDistFlow} and substitute them into the DistFlow equations \eqref{DistFlow}.  We then solve the DistFlow equations exactly and compare the relative losses (normalized with respect to the losses of the no-optimization scheme) of the DistFlow (Loss$^{\rm{glob}}_{\rm DF}$) and LinDistFlow equations (Loss$^{\rm glob}_{\rm lin}$).
We also compare the maximum difference between the normalized voltages calculated by LinDistFlow (V$^{\rm glob}_{\rm lin}$) and by DistFlow (V$^{\rm{glob}}$) equations for each case. We use the Matpower package \cite{matpower} to solve the DistFlow solutions, with the reactive powers set to values obtained by solution of the LinDistFlow optimization problem \eqref{LinDistFlow}.  The table below shows these comparisons, with $\left| \delta V \right|$ defined as $\max_j \left| {V^{\rm glob}_{\rm lin}}_j - V^{\rm{glob}}_j  \right|$.

In addition, we also solve the global optimization DistFlow problem with CVX package using a centralized algorithm from \cite{farivar2012optimal}, which finds the values of reactive power $q^{(g)}$ that solve Eqs.~\eqref{DistFlowOptimization}. The obtained relative losses for such a solution (shown in the table below as Loss$^{\rm{glob-nonlin}}$) are close to the relative losses calculated for the reactive powers $\tilde{q}_j^{(g)}$ obtained from our decentralized optimization of LinDistFlow equations.

%For the cases 6 and 7, where the voltage constraints are violated, the
%actual minimum (maximum) voltages are slightly higher than the ones calculated from the LinDistFlow equations. Hence the relative losses Loss$^{\rm glob}_{\rm{DF}}$ are slightly lower for case 6,

\begin{center}
\begin{tabular}{|c|cccc|}
\hline
 {\bf Case} & {\bf Loss$^{\rm glob}_{\rm lin}$} & {\bf Loss$^{\rm glob}_{\rm{DF}}$} & {\bf Loss$^{\rm glob-nonlin}$}& $\left| \delta V \right|$ \\
 \hline
 1 &    $0.834$  &  $0.828$ &  $0.828$  & $0.000$  \\
 2 &    $0.941$  &  $0.937$ &  $0.937$ & $0.000$  \\
 3 &    $0.847$ &   $0.821$ &  $0.819$  & $0.001$ \\
 4 &    $0.948$   & $0.937$ &  $0.937$ & $0.000$ \\
 5 &   $0.700$ &    $0.687$ & $0.687$   & $0.000$ \\
 6 &   $0.954$ &    $0.916$ & $0.920$  & $0.003$ \\
 7 &   $1.087$ &    $1.112$ & $1.071$& $0.005$ \\ \hline

\end{tabular}
\end{center}

The relative losses Loss$^{\rm glob}_{\rm{DF}}$ agree within few percent with the values predicted by the LinDistFlow equations, and in most cases, the DistFlow relative losses are actually slightly lower than the ones obtained from LinDistFlow equations. For all the cases considered, the maximum difference between values of voltages are of the order of $10^{-3}$, which is a satisfactory precision since the algorithm needs to optimize objective function subject to (normalized) voltage constrained to between $0.95$ and $1.05$.

For the cases 6 and 7, where global reactive power control scheme is necessary to ensure voltage regulation within allowed bounds, the voltage magnitude calculated with LinDistFlow equations ($V^{\rm glob}_{\rm lin}$) slightly overestimates the exact value calculated from DistFlow equations ($V^{\rm glob}$). As a result, the actual voltage is in fact lower by $O(10^{-3})$. Hence, the global optimization of LinDistFlow equations leads to lower relative losses (Loss$^{\rm glob}_{\rm{DF}}$) than Loss$^{\rm glob-nonlin}$ obtained from optimization of DistFlow problem \eqref{DistFlowOptimization} for case 6 (where LinDistFlow optimization's minimum voltage is $0.947$), and higher relative losses for case 7, as the maximum normalized voltage obtained from LinDistFlow optimization is actually $1.045$, lower than the maximum bound $1.05$, which is attained with the global DistFlow optimization.

\section{Discussion}
\label{sec_disc}
In this work, we presented three distributed algorithms achieving global optimality of reactive power flows in power-distribution systems and tested these algorithms on multiple examples of circuits with high PV penetration. In this formulation of the  optimization, the cost function and constraints are separable and this structure is explored to construct and compare three exact distributed algorithms that rely on local measurements, local computations and communications between nearest neighbors. The general advantage of this distributed implementation is that each node can have its own set of constraints as well as its own objective function, which do not need to be known to the other nodes. The dual-ascent algorithm allows for simple implementation, but in most cases requires significantly more iterations to converge. The consensus-based ADMM algorithms require more sophisticated  local computations, however, overall they converge much faster than the dual ascent algorithm.

Our approach allows multiple generalizations.  We highlight some, while emphasizing the importance of future theoretical, algorithmic and experimental (simulations and testbed) explorations:
\begin{itemize}
\item One can easily generalize our ADMM algorithm to different objective functions: for instance, some nodes can have an additional term in the objective function which depends on voltage, such as $\left( V - V_0 \right)^2$. This correction would improve the voltage quality for that particular node. Those rules can even be changed locally and on-the-go, without other nodes being aware of the adaptation.

\item Communications between different nodes do not need to be synchronized and may allow some degree of delay.  It will be important to analyze the robustness of the scheme to delays, errors and corruptions, e.g. originating from targeted attacks.  It is expected that the distributed nature of the algorithm makes it much more secure against localized and non-correlated attacks than centralized algorithms.

\item Our algorithm is naturally extendable to account for more complicated tree-like topologies and some (still separable) nonlinearities in power flows, of the type already discussed in \cite{Zampieri13c} for the dual ascent algorithm. It also allows generalization to account for binary PV-inverter selection such as discussed in \cite{Dall2013}.

\end{itemize}

\section{Acknowledgments}
The authors thank Sidhant Misra, Changhong Zhao and Steven Low for helpful advice and productive discussions. The work at LANL was funded by the Advanced Grid Modeling Program in the Office of Electricity in the US Department of Energy and was carried out under the auspices of the National Nuclear Security Administration of the U.S. Department of Energy at Los Alamos National Laboratory under Contract No. DE-AC52-06NA25396. The work at NMC is supported by the National Science Foundation award \# 1128501, EECS Collaborative Research ``Power Grid Spectroscopy".

\bibliographystyle{ieeetr}
\bibliography{distcontrol}

\end{document}